\documentclass[12pt, oneside]{article}
\usepackage{latexsym}
\usepackage{graphicx}
\usepackage{setspace}
\usepackage[all]{xy}
\usepackage{geometry}
\usepackage[psamsfonts]{amsfonts}
\usepackage{amsmath, amsthm, amssymb}
\usepackage{epsfig}
\usepackage{makeidx}
\usepackage{enumerate}
\usepackage{color}

\xyoption{arc}

\usepackage{verbatim}

\newcommand{\IZ}{\mathbb Z}

\newcommand{\IC}{\mathbb C}

\newtheorem{Thm}{Theorem}
\newtheorem{Cor}{Corollary}
\newtheorem{Def}{Definition}
\newtheorem{Prop}{Proposition}

\newtheorem{Lemma}{Lemma}
\newtheorem{Notation}{Notation}

\theoremstyle{definition}
\newtheorem{Remark}{Remark}
\newtheorem{Example}{Example}

\advance \topmargin by -\headheight
\advance \topmargin by -\headsep
\evensidemargin \oddsidemargin
\begin{document}

%% TITLE

\title{Some graph properties determined by edge zeta functions}
\author{Christopher Storm\\
  Department of Mathematics and Computer Science,\\
  Adelphi University,\\
  \texttt{cstorm@adelphi.edu}}
\date{\today}
\maketitle

%% ABSTRACT

\begin{center}{\bf \large Abstract}\end{center}
\begin{quote}
Stark and Terras introduced the edge zeta function of a finite graph in 1996.  The edge zeta function is the reciprocal of a polynomial in twice as many variables as edges in the graph and can be computed in polynomial time.  We look at graph properties which we can determine using the edge zeta function.  In particular, the edge zeta function is enough to deduce the clique number, the number of Hamiltonian cycles, and whether a graph is perfect or chordal.  Actually computing these properties takes exponential time.  Finally, we present a new example illustrating that the Ihara zeta function cannot necessarily do the same.\\
\end{quote}

\section{Introduction}

In 1996, Stark and Terras introduced the edge zeta function of a finite graph as a generalization of the Ihara zeta function \cite{MR1399606}.  Horton, Stark, and Terras used the edge zeta function in 2006 \cite{MR2277616} to provide a new proof of Bass's determinant expression for the Ihara zeta function \cite{MR1194071}.  Aside from this, the edge zeta function hasn't received much attention.  Our goal is to show that the edge zeta function, which can be computed in polynomial time, determines a large amount of information about a graph.  We hope that this can then be used to show that this invariant is very good at distinguishing graphs.\\

For the rest of this section, we give the definition of the edge zeta function and the Ihara zeta function.  In Section 2, we survey some known properties of graphs which are determined by the Ihara zeta function.  Then in Section 3, we look specifically at the edge zeta function.  We will show that the edge zeta function determines the clique number, the number of Hamiltonian cycles, and the presence or absence of holes and antiholes in a graph, allowing us to conclude if a graph is perfect or chordal.\\

We begin by defining graphs, digraphs, and the symmetric digraph associated to a graph.  All structures treated here are finite.   We refer the reader to the books by Harary, and Chartrand and Lesniak \cite{MR0256911, MR834583} for a good overview of these structures.\\

A \emph{graph} $X = (V, E)$ is a finite nonempty set $V$ of \emph{vertices} and a finite multiset $E$ of unordered pairs of vertices, called \emph{edges}.  If $\{u, v\} \in E$, we say that $u$ is \emph{adjacent} to $v$ and write $u \sim v$.  A graph $X$ is \emph{simple} if there are no edges of the form $\{v, v\}$ and if there are no repeated edges.\\

A \emph{directed graph} or \emph{digraph} $D = (V, E)$ is a finite nonempty set $V$ of \emph{vertices} and a finite multiset $E$ of ordered pairs of vertices called \emph{arcs}.  For an arc $e = (u, w)$, we define the \emph{origin} of $e$ to be $o(e) = u$ and the \emph{terminus} of $e$ to be $t(e) = w$.  The \emph{inverse arc} of $e$, written $\overline{e}$, is the arc formed by switching the origin and terminus of $e$: $\overline{e} = (w, u)$.  In general, the inverse arc of an arc need not be present in the arc set of a digraph.\\

\begin{figure}
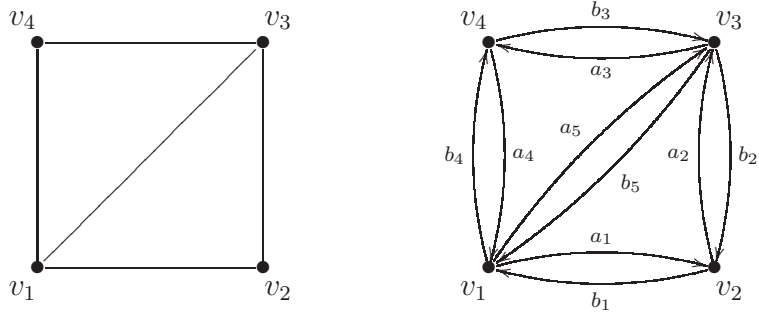

\[\xy
(0,0)*!C++\xybox{
(0,0)*{\bullet} = "v1"+(-2, -3)*{v_1};
(30, 0)*{\bullet} = "v2"+(2, -3)*{v_2};
(0, 30)*{\bullet} = "v3"+(-2,3)*{v_4};
(30,30)*{\bullet} = "v4"+(2, 3)*{v_3};
{\ar@{-} "v1"; "v2"};
{\ar@{-} "v2"; "v4"};
{\ar@{-} "v4"; "v3"};
{\ar@{-} "v3"; "v1"};
{\ar@{-} "v1"; "v4"};
};
(60,0)*!C++\xybox{
(0,0)*{\bullet} = "v1"+(-2, -3)*{v_1};
(30, 0)*{\bullet} = "v2"+(2, -3)*{v_2};
(0, 30)*{\bullet} = "v3"+(-2,3)*{v_4};
(30,30)*{\bullet} = "v4"+(2, 3)*{v_3};
{\ar^{a_1}@{->}@/^/  "v1"; "v2"};
{\ar^{a_2}@{->}@/^/  "v2"; "v4"};
{\ar^{a_3}@{->}@/^/  "v4"; "v3"};
{\ar^{a_4}@{->}@/^/  "v3"; "v1"};
{\ar^{a_5}@{->}@/^/  "v1"; "v4"};
{\ar^{b_1}@{->}@/^/  "v2"; "v1"};
{\ar^{b_2}@{->}@/^/  "v4"; "v2"};
{\ar^{b_3}@{->}@/^/  "v3"; "v4"};
{\ar^{b_4}@{->}@/^/  "v1"; "v3"};
{\ar^{b_5}@{->}@/^/  "v4"; "v1"};
};
\endxy\]
\caption{The complete graph minus an edge and its symmetric digraph}
\label{fig:SymGraph}
\end{figure}

A digraph $D$ is called \emph{symmetric} if, whenever $(u,w)$ is an arc of $D$, its inverse arc $(w, u)$ is as well.  There is a natural one-to-one correspondence between the set of symmetric digraphs and the set of graphs, given by identifying an edge of the graph to an arc and its inverse arc on the same vertices.  We denote by $D(X)$ the symmetric digraph associated with the graph $X$.  We give an example in Figure \ref{fig:SymGraph}.\\

To define the zeta functions, we need several cycle definitions.  We let $X$ be a graph and $D(X)$ its symmetric digraph.  A \emph{cycle $c$ of length $n$} in $X$ is a sequence $c = (e_1, \cdots, e_n)$ of $n$ arcs in $D(X)$ such that $t(e_i) = o(e_{i+1})$ for $1 \leq i \leq n-1$ and $t(e_n) = o(e_1)$.  We say that $c$ has \emph{backtracking} if $\overline{e_{i+1}} = e_i$ for some $i$ satisfying $1 \leq i \leq n-1$.  Also, $c$ has a \emph{tail} if $e_1 = \overline{e_n}$.  We are primarily interested in cycles with no backtracking or tails.\\

The \emph{$r$-multiple} of the cycle $c$ is the cycle $c^r$ formed by going $r$ times around $c$.  We say a cycle is \emph{primitive} if it is not the $r$-multiple of some other cycle $b$ for $r \geq 2$.  We impose an equivalence relation on cycles via cyclic permutation; i.e., two cycles $b = (e_1, \cdots, e_n)$ and $c = (f_1, \cdots, f_n)$ are \emph{equivalent} if there is a fixed $\alpha \in \IZ / n\IZ$ such that $e_i = f_{i + \alpha}$ for all $i \in \IZ / n\IZ$ (all indices are considered mod $n$).  Note that the direction of travel does matter so that traversing a cycle in the opposite direction is not equivalent to the original cycle.  A \emph{prime cycle} is the equivalence class of primitive cycles which have no backtracking or trails, written $[c]$.\\

For a graph $X$ with symmetric digraph $D(X)$, we associate to each arc $e$ of $D(X)$ an invariant $u_e$.  Then for a prime cycle $[c]$, we define a function
\[g(c) = \prod_{\text{$e$ arc in $c$}} u_e.\]
This function reports which arcs are used in a prime cycle and how many times they are used.\\
\begin{Example}
Use the labeling given in Figure \ref{fig:SymGraph}.  Then the cycles described by $\{a_1, a_2, b_5, a_1, a_2, b_5, b_4, b_3, b_5, b_4, b_3, b_5\}$ and $\{a_1, a_2, b_5, b_4, b_3, b_5, b_4, b_3, b_5 ,a_1, a_2, b_5\}$ both have
\[g(c) = u_{a_1}^2 u_{a_2}^2 u_{b_5}^4 u_{b_4}^2 u_{b_3}^2.\]
\end{Example}

We can now define the edge and Ihara zeta functions of a graph:\\
\begin{Def}[Stark and Terras]
For a finite graph $X$, associate to each arc of $D(X)$ an invariant $u_e$.  The \emph{edge zeta function} of $X$ is a function of $u_e \in \IC$ (sufficiently near $0$) given by
\[\zeta_X(\vec{u}) = \prod_{\text{primes cycles $[c]$}}\left(1-g(c)\right)^{-1}.\]
The \emph{Ihara zeta function} of $X$ is given by specializing each $u_e$ to $u$, which is
\[Z_X(u) = \prod_{\text{primes cycles $[c]$}}\left(1 - u^{l(c)}\right)^{-1},\]
where $l(c)$ is the \emph{length} of a representative of the prime cycle $[c]$.\\
\label{Def:EdgeZeta}
\end{Def}

Remarkably, the edge zeta function of a finite graph is the reciprocal of a multivariate polynomial.  To see this, we define the \emph{directed edge matrix} $T$ associated to a graph.  For a graph $X$, we begin by fixing a labeling of the arcs of $D(X)$.\\

\begin{Def}[Stark and Terras]
The \emph{directed edge matrix} $T$ has as its $ij$ entry
\[t_{ij} = 
\begin{cases}
1& \text{if $t(e_i) = o(e_j)$ and $e_i \neq \bar{e_j}$;} \\
0 & \text{otherwise.}
\end{cases}\]
We let $U$ be the diagonal matrix containing the indeterminants from Definition \ref{Def:EdgeZeta}:
\[U = \rm{diag}(u_{e_1}, \cdots, u_{e_{2|E|}}).\]
\label{Def:EdgeMat}
\end{Def}

We note that other authors have relied upon this $T$ matrix as well.  Kotani and Sunada \cite{MR1749978} use it as the Perron--Frobenius operator of the oriented line graph associated to $X$.  From the matrices in Definition \ref{Def:EdgeMat}, we realize a determinant expression for the edge zeta function (and thus for the Ihara zeta function as well).\\
\begin{Thm}[Stark and Terras]
Let $X$ be a finite graph.  With the notation of Definitions \ref{Def:EdgeZeta} and \ref{Def:EdgeMat}, we have
\[\zeta_X(\vec{u})^{-1} = \det(I - UT) = \det(I - TU).\]
\label{Thm:FirstDet}
\end{Thm}

Hence the edge zeta function is the reciprocal of a multivariate polynomial in at most $2|E(X)|$ variables and can be computed in polynomial time.  In addition, and very importantly for us, given the edge zeta function of a graph $X$, it is very easy to specialize it to realize the edge zeta function of a subgraph of $X$.\\

\begin{Prop}[Stark and Terras]
Let $X$ be a graph with symmetric digraph $D(X)$.  Let $F$ be a subset of $E(X)$, and let $\mathbb{F}$ consist of the set of arcs in $D(X)$ corresponding to the edges in $F$.  Suppose $W$ is the graph obtained from $X$ by erasing all of the edges in $F$.  Then
\[\zeta_X(\vec{u})|_{u_e=0, \forall e \in F} = \zeta_W(\vec{u}).\]
\label{Prop:EdgeDeletions}
\end{Prop}

We will use Proposition \ref{Prop:EdgeDeletions} over and over again in the final section.  It will be our main tool for picking out graph properties based on the edge zeta function.  Our general technique is to identify graphs which are uniquely determined by their Ihara zeta function.  Then, with the aid of Proposition  \ref{Prop:EdgeDeletions} we can count how many subgraphs have the desired zeta function---and are thus determined.  We will assume throughout that we are given the identification of the indeterminant of an arc and its inverse arc, so that we can directly specialize to get edge zeta functions of subgraphs.\\

In the next section, we survey the properties of graphs that are known to be determined by the Ihara zeta function.  Then, in Section 3, we will look specifically at edge zeta functions and see how we can realize more graph invariants.\\

% PROPERTIES FROM IHARA ZETA
\section{Properties determined by the Ihara zeta function}

In this section, we look at some of the known results about the single variable Ihara zeta function which will prove useful to us in Section 3.  We begin by exploring some of the consequences of Theorem  \ref{Thm:FirstDet}.  Then, we look at a more detailed determinant expression, given by Bass, and see that regular graphs are \emph{cospectral} if and only if they have the same zeta function.  This last fact will be very useful at identifying structure determined by the edge zeta function.\\

We now take a closer look at Theorem \ref{Thm:FirstDet}.  For a graph $X$, the Ihara zeta function $Z_X(u)$ can be written as $\det(I - uT)^{-1}$ where $T$ is the directed edge matrix associated with $X$.  From this expression, one can deduce that the maximum degree of the reciprocal of the zeta function is $2|E(X)|$.  In fact, if there are no vertices of degree 1 in $X$, this is exactly the degree of the polynomial.  This fact has been noted by Stark and Terras \cite{MR1399606} as well as by Kotani and Sunada \cite{MR1749978}.  Czarneski gives a proof of this by computing the $2|E(X)|^{th}$ coefficient of the reciprocal of the zeta function and showing that it's non-zero so long as all of the vertices have degree at least 2.  In addition, Horton gives a detailed discussion of the eigenvalues and eigenfunctions of the matrix $T$ which sheds light on this fact \cite{Horton-2006}.  Hence, for a finite graph where every vertex has at least degree 2---such a graph will be refered to as \emph{md2} from now on---the zeta function determines the number of edges in the graph.\\

What happens when a vertex has degree 1?  Recalling our prime cycle definitions given in the previous section, the only way to include an edge which is incident to a degree 1 vertex in a cycle is to either have backtracking or a tail.  Hence, any edges incident to a degree 1 vertex are completely ignored by the zeta function.  One can then remove these edges and vertices.  This may create new edges of degree 1, which can also be removed, successively, until the remaining graph is $md2$.  This underlying graph is what the zeta function is really studying.\\

\begin{Notation}
For a graph $X$, we let $m = |E|$.  We write
\[\frac{1}{Z_X(u)} = Z_X(u)^{-1} = c_0 + c_1u + c_2u^2 + c_3u^3 + \cdots + c_{2m}u^{2m}.\]
We denote by $c_k(X)$ the coefficient $c_k$ of $u^k$ of the reciprocal of $Z_X(u)$.\\
\label{Not:ck}
\end{Notation}

For a detailed discussion of how the numbers $c_k(X)$ relate to the structure of $X$, we refer the reader to \cite{scott-2007}.  We will be particularly interested in $c_{2m}(X)$.  This coefficient will depend on the degree sequence of $X$ as detailed, independently, by Czarneski \cite{Czarneski-2005} and Horton \cite{Horton-2006}.\\

\begin{Prop}[Czarneski, 2005; Horton, 2006]
Let $X$ be a finite graph with \\$|E(X)| = m$.  Then
\[c_{2m}(X) = \prod_{v \, \in \, V(X)}\left( d(v) - 1\right).\]
\label{Prop:BigCoeff}
\end{Prop}

We now look at a more detailed determinant expression, given by Bass, which generalizes Ihara's initial determinant expression \cite{MR0223463} of the zeta function of a regular graph.\\

\begin{Thm}[Bass]
Let $X$ be a finite, connected graph with adjacency matrix $A$ and degree matrix $D$ defined as a diagonal matrix with the degrees of the vertices of $X$ down the diagonal.  Let $I$ be the identity matrix.  Then,
\[Z_{X}(u) = (1 - u^2)^{\chi(X)} \det(I - uA + u^2(D-I))^{-1}\]
where $\chi(X) = |V| - |E|$ is the Euler Number of the graph $X$.\\
\label{Thm:Bass}
\end{Thm}

A great deal of the theory of Ihara zeta functions comes from a study of Theorem \ref{Thm:Bass}.  The following observation---first made by Mellein \cite{Mellein} although certainly known to Quenell \cite{Quenell-1998}---is very useful for us.\\

\begin{Thm}[Mellein]
Suppose $X$ and $Y$ are both $k$-regular graphs.  Then $X$ and $Y$ are \emph{cospectral}---their adjacency matrices have the same spectra---if and only if
\[Z_X(u) = Z_Y(u).\]
\label{Thm:Mell}
\end{Thm}
\begin{proof}
We only give a very broad sketch of the proof.  Since $X$ and $Y$ are both $k$-regular, the number of vertices in $X$ and $Y$ can be determined based on the number of edges.  Then we only need to study the determinant expression that appears in Theorem \ref{Thm:Bass}.\\

Since the graphs are regular, all of the matrices inside the determinant commute.  This allows us to simultaneously diagonalize, giving us
\[\det(I - uA + qu^2I) = \prod_{\lambda_i \in {\rm spec }A}(1 - \lambda_iu + qu^2),\]
where $q = k-1$ and $A$ is the adjacency matrix of $X$ or $Y$ as needed.  Manipulating this last expression gives us the result.\\
\end{proof}

Theorem \ref{Thm:Mell} will be a key fixture in Section 3.  Whenever a $k$-regular graph is uniquely determined by its spectrum, we will be able to conclude that its Ihara zeta function is also uniquely determined.  This will allow us to search for specific structures which could appear as subgraphs in a graph.\\

\begin{figure}[th]
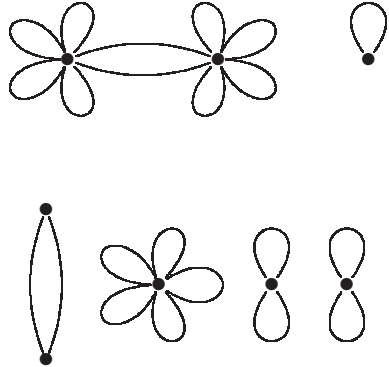

\[\xy
(0,0)*!C++\xybox{
(0,0)*{\bullet} = "v1";
(20, 0)*{\bullet} = "v2";
(40, 0)*{\bullet} = "v3";
{\ar@{-}@/^/ "v1"; "v2"};
{\ar@{-}@/^/ "v2"; "v1"};
"v1"; "v1" **\crv{ (-7,0) & (-10, 7) & (-2,5) };
"v1"; "v1" **\crv{ (-7,0) & (-10, -7) & (-2,-5) };
"v1"; "v1" **\crv{ (-2,5) & (2, 10) & (5,5) };
"v1"; "v1" **\crv{ (-2,-5) & (2, -10) & (5,-5) };
"v2"; "v2" **\crv{ (22,5) & (18, 10) & (15,5) };
"v2"; "v2" **\crv{ (22,-5) & (18, -10) & (15,-5) };
"v2"; "v2" **\crv{ (27,0) & (30, 7) & (22,5) };
"v2"; "v2" **\crv{ (27,0) & (30, -7) & (22,-5) };
"v3"; "v3" **\crv{ (36,5) & (40, 10) & (44,5) };
};
(0,-30)*!C++\xybox{
(0,0)*{\bullet} = "v1";
(0, -20)*{\bullet} = "v2";
(15, -10)*{\bullet} = "v3";
(30,-10)*{\bullet} = "v4";
(40,-10)*{\bullet} = "v5";
{\ar@{-}@/^/ "v1"; "v2"};
{\ar@{-}@/^/ "v2"; "v1"};
"v5"; "v5" **\crv{ (36,-5) & (40, 0) & (44,-5) };
"v5"; "v5" **\crv{ (36,-15) & (40, -20) & (44,-15) };
"v4"; "v4" **\crv{ (26,-5) & (30, 0) & (34,-5) };
"v4"; "v4" **\crv{ (26,-15) & (30, -20) & (34,-15) };
"v3"; "v3" **\crv{ (8,-10) & (5, -3) & (13,-5) };
"v3"; "v3" **\crv{ (8,-10) & (5, -17) & (13,-15) };
"v3"; "v3" **\crv{ (13,-5) & (17, 0) & (20,-5) };
"v3"; "v3" **\crv{ (13,-15) & (17, -20) & (20,-15) };
"v3"; "v3" **\crv{ (20,-6) & (27, -10) & (20,-14) };
};
\endxy\]
\caption{Two graphs with the same zeta function but different numbers of vertices and connected components.}
\label{fig:CzarExe}
\end{figure}

\begin{Example}
Lest we get too carried away, we give an example as a warning.  Czarneski \cite{Czarneski-2005} gave an example of a pair of graphs which have the same zeta function but differing numbers of vertices and connected components.  Cooper \cite{Yaim-2006} has extended this example to give an infinite family of pairs of such graphs.  Czarneski's original example is found in Figure \ref{fig:CzarExe}.\\
\end{Example}

The results given above are all that we will need in Section 3.  We would be remiss in not mentioning that this section is not exhaustive.  In particular, the Ihara zeta function determines the number of spanning trees in a graph \cite{MR1080105, MR1654160} in of an analogue to the class number formula of a number field.  It also determines the \emph{girth} of a graph \cite{Horton-2006, scott-2007}\\

In the event that a graph is $k$-regular, we will be able to determine whether or not it is connected from its edge zeta function.  We begin the next section by showing how to determine whether or not a graph is regular, and if it is, how to determine connectivity.  Once these preliminaries are out of the way, we go straight into counting subgraphs in the graph.\\

%EDGE ZETA PROPERTIES
\section{Properties determined by the edge zeta function}

We saw in the previous section that we cannot necessarily count the number of vertices or connected components of a graph $X$ just from its Ihara zeta function.  We begin this section by showing how to determine whether or not a graph is regular from its edge zeta function.  Once we've concluded that a graph is regular, we will be able to determine whether or not it is connected.  From this point, we will be able to make use of Theorem \ref{Thm:Mell} in conjunction with results about regular graphs which are determined by their adjacency matrix spectra to identify the properties we desire.\\

%Determining regularity
\begin{Lemma}
Let $X$ be a finite md2 graph with $|E(X)| = m$.  Suppose there exists a vertex $v \in V(X)$, satisfying $d(v) > 2$, with $v$ adjacent to two vertices $x, y \in V(X)$ such that
\[d(x) \neq d(y).\]
Then 
\[c_{2m-2}(X \setminus \{v, x\}) \neq c_{2m-2}(X \setminus \{v, y\}),\]
where $X \setminus e$ is the graph formed by removing edge $e$.\\
\label{Lem:3degs}
\end{Lemma}
\begin{proof}
This follows directly from the definitions in Notation \ref{Not:ck} and Proposition \ref{Prop:BigCoeff}.\\
\end{proof}

We now show how to determine whether a graph is regular, biregular bipartite, bipartite with all of the vertices in one set having degree 2, or none of the above.\\  

\begin{Prop}
Let $W$ be a md2 graph with $m$ edges.  We denote by $W \setminus e$ the subgraph of $W$ which is formed by removing the edge $e$.  Suppose that the numbers
\[\{c_{2m-2}(W \setminus e_1), c_{2m-2}(W \setminus e_2), \cdots, c_{2m-2}(W \setminus e_m)\}\]
are all the same.  Then $W$ satisfies one of the following:
\begin{enumerate}
\item $W$ is $k$-regular, and $k$ can be determined.
\item $W$ is a $(p,q)$-biregular bipartite graph for some $p,q \geq 3$.
\item $W$ is a bipartite graph where all of the vertices in one of the vertex sets have degree 2.
\end{enumerate}
\label{Prop:Regs}
\end{Prop}
\begin{proof}

We first note that $W$ is $2$-regular if and only if $c_{2m}(W) = 1$.\\

We assume for the moment that $W$ is connected and that $d(v) \geq 3$ for all $v \in V(W)$.  Suppose that two vertices $v_1$ and $v_2$ are adjacent and satisfy $d(v_1) = d(v_2) = k$; then, we claim $W$ is $k$-regular.  For any other vertex $w$, we consider a path from $v_1$ to $w$ given by $\{v_1 = w_1, w_2, \cdots, w_n = w\}.$  Then $v_1$ is adjacent to $v_2$ and $w_2$ (it's ok if they are the same vertex), so $d(w_2) = d(v_2)$, or we fail the conditions of the proposition because of Lemma \ref{Lem:3degs}.  Similarly, $d(w_3) = d(w_1)$, and we continue until $d(w) = d(w_{n-2}) = d(v_1)$.  Thus $W$ is $k$-regular.\\

Now suppose there are two vertices which are adjacent and satisfy $d(v_1) = p$ and $d(v_2) = q$ with $p \neq q$.  By a similar argument, we see that every vertex must have degree $p$ or degree $q$.  Now suppose that there exists an odd cycle in $W$.  Either two vertices of the same degree are adjacent --- forcing the graph to be $k$-regular: a contradiction --- or there is a third degree: another contradiction.  Thus every cycle must have even length, and $W$ is in fact $(p,q)$-biregular bipartite.\\

We can remove the condition on connectivity.  If there are more than one connected component, then each component must have the same degree structures.  Else removing an edge in one component and an edge in a different component would give different numbers $c_{2m-2}$ for those edge removals.\\

We distinguish between these two cases.  If $W$ is $k$-regular, then 
\[c_{2m-2}(W \setminus e) = \left[\frac{k-2}{k-1}\right]^2 c_{2m}(W),\]
for all $e \in E(W)$.  If $W$ is $(p,q)$-biregular bipartite, then
\[c_{2m-2}(W \setminus e) = \left[\frac{p-2}{p-1}\right] \left[\frac{q-2}{q-1}\right] c_{2m}(W),\]
for all $e \in E(W)$.  From these expressions we can distinguish which case we have, and if the graph is regular determine the value $k$.\\

These cases cover the situation when $c_{2m-2}(W \setminus e)$ is non-zero.  If $c_{2m-2}(W \setminus e) = 0$ for all $e \in E(W)$; then, every edge must be incident to a vertex of degree 2.  In this case, we can identify $W$ as belonging to category 3 above.
\end{proof}

%Determining connected in regular graphs.  
Now that we can distinguish whether a graph is $k$-regular or not, we show how to tell when a $k$-regular graph is connected.\\

\begin{Prop}[Connectivity in regular graphs]
Suppose $X$ is a $k$-regular graph.  Then $X$ is connected if and only if the pole of $Z_X(u)$ at $u = \frac{1}{k-1}$ is simple.\\
\label{Prop:Connect}
\end{Prop}
\begin{proof}
The multiplicity of the pole of $Z_X(u)$ at $u = \frac{1}{k-1}$ is 1 if and only if the multiplicity of $\lambda = k$ as an eigenvalue of the adjacency matrix of $X$ is 1.  This is true if and only if $X$ is connected.\\
\end{proof}

We are now ready to look at some properties of graphs which are determined by the edge zeta function.  Our method is simple.  We will use the previous propositions in conjunction with Theorem \ref{Thm:Mell} to identify edge-induced and vertex-induced subgraphs which are isomorphic to particular graphs.  We first establish some useful notation.\\

\begin{Def}
Let $X = (V, E)$ be a finite graph.  For a subset $S$ of $V$, the \emph{vertex-induced subgraph $\langle S \rangle$} of $X$ is the subgraph formed by taking $S$ as its vertex set and taking the set of edges which have both endpoints in $S$ as the edge set.  For a subset $R$ of $E$, the \emph{edge-induced subgraph $\langle R \rangle$} of $X$ is the subgraph formed by taking $R$ as its edge set and the set of vertices which are incident to some edge in $R$ as the vertex set.\\

For a graph $W$, we denote by $s_v(W, X)$ the number of vertex-induced subgraphs of $X$ which are isomorphic to $W$.  Similarly, we denote by $s_e(W, X)$ the number of edge-induced subgraphs of $X$ which are isomorphic to $W$.\\
\label{Def:InducedCounts}
\end{Def}

We now give our main theorem, which will drive the rest of the section.\\

\begin{Thm}[Counting subgraphs]
Let $X$ be a md2 graph with edge zeta function $\zeta_X(\vec{u})$.  Let $W$ be a $k$-regular graph which is determined by the spectrum of its adjacency matrix.  Then the numbers $s_v(W, X)$ and $s_e(W,X)$ are both determined by $\zeta_X(\vec{u})$.\\
\label{Thm:MainThm}
\end{Thm}

\begin{proof}
We first show how to determine $s_e(W, X)$.  Suppose that $|E(W)| = \tilde{m}$.  We denote by $\mathcal{S}$ the set consisting of all unordered $\tilde{m}$-tuples of the arc/inverse arc pairs of indeterminants that appear in $\zeta_X$.\\

For an element $R \in \mathcal{S}$ we form the function $\zeta_{\langle R \rangle}$.  Due to Proposition \ref{Prop:EdgeDeletions}, this is exactly the zeta function of the edge-induced subgraph of $X$ given by the edges indexed in $R$.  We use Proposition \ref{Prop:Regs} to verify that the edge-induced subgraph is a regular graph.  We specialize to its Ihara zeta function and then use Theorem \ref{Thm:Mell} to check if the edge-induced subgraph is isomorphic to $W$ or not.  We repeat this process for every element of $\mathcal{S}$ to compute $s_e(W, X)$.\\

Interestingly, with full use of the edge zeta function, it is not much more difficult to compute $s_v(W, X)$.  Suppose $R$ is a subset of $\mathcal{S}$ which contributed to $s_e(W, X)$.  We now pick an edge $e$ which isn't represented in $R$.  Then, we form the edge zeta function induced from the set $R\cup\{e\}$.  We specialize to the Ihara zeta function of this graph.  Now, there are three options for how $e$ interacts with the edge induced subgraph of $X$ which comes from $R$.  If $e$ is incident to zero or one vertices incident to an edge in $R$, the Ihara zeta function will be exactly the Ihara zeta function that arose just from $R$.  If, however, $e$ is incident to two vertices which are incident to edges in $R$, the Ihara zeta function will change.  In particular, its maximum degree will increase by $2$.\\

To compute $s_v(W, X)$, we simply pick each subset $R$ of $\mathcal{S}$ and then perform the above process with each edge not in $R$.  If the Ihara zeta function of the new graphs always matches the one induced from $R$, we have a vertex-induced subgraph isomorphic to $W$.  If it does change for any edge, we don't.\\  
\end{proof}

\begin{Remark}
A slightly more general statement of Theorem \ref{Thm:MainThm} is possible.  Czarneski \cite{Czarneski-2005} gave a statement of Theorem \ref{Thm:Mell} to biregular bipartite graphs.  Using this statement, we could also consider graphs $W$ which are biregular bipartite and uniquely determined by the spectrum of their adjacency matrix.\\
\end{Remark}

Theorem \ref{Thm:MainThm} provides a machine to identify substructures in $X$.  The study of graphs which are determined by their spectra is an old one, dating back to chemistry in 1956 \cite{Primas-1956}.  Fisher \cite{MR0197344} also addressed this question in response to Kac's \cite{MR0201237} famous question ``Can one hear the shape of a drum?"  We recommend the excellent book by Biggs \cite{MR1271140} and article by van Dam and Haemers \cite{MR2022290} as a starting point to the literature on these questions.\\

We will focus on complete graphs and cycles as they play important roles in determining the structure of a graph.  The \emph{complement of a graph} $X$ is the graph $\bar{X}$ formed by keeping the same vertex set and edge set formed by making $\{u,v\}$ an edge in $\bar{X}$ whenever it is not one in $X$.  The following proposition, which can be found in \cite{MR2022290}, is straight-forward.\\

\begin{Prop}
The complete graph $K_n$, the cycle $C_n$ and their complements are determined by their adjacency matrix spectrum.
\label{Prop:DSGraphs}
\end{Prop}

We look at the graphs in Proposition \ref{Prop:DSGraphs} individually.  We say that a graph $X$ with $n$ vertices is \emph{Hamiltonian} if it has an edge-induced subgraph isomorphic to $C_n$.  Such a cycle is called a Hamiltonian cycle, and we denote by Ham$(X)$ the number of such cycles in $X$.\\

%C_n is determined.
%Hamiltonian cycles are determined
\begin{Cor}[Cycles and Hamiltonian cycles]
Let $X$ be a graph on $n$ vertices.  Then for $k = 3, \cdots, n$, the number $s_e(C_k, X)$ is determined by $\zeta_X(\vec{u})$.  In particular, $\rm{Ham}(X)$ is determined.\\
\end{Cor}

In fact, the numbers $s_v(C_k, X)$ will also be very interesting.  We return to these numbers in a moment, after we look at counting copies of complete graphs in $X$.  The \emph{clique number} of a graph $X$, written $\omega(X)$ is the largest integer $r$ such that $X$ has a vertex-induced subgraph isomorphic to $K_r$.  The clique number is often associated with coloring as it gives an immediate lower bound on the chromatic number.\\

\begin{Cor}[Complete graphs and the clique number]
Let $X$ be a graph on $n$ vertices.  Then for $r=3, \cdots, n$, the number $s_e(K_r, X) = s_v(K_r, X)$ is determined by $\zeta_X(\vec{u})$.  In particular, $\omega(X)$ is determined.\\
\end{Cor}

We mention two important classes of graphs since their structure is dependent upon the presence or absence of copies of $C_k$ and $\bar{C_k}$ as vertex-induced subgraphs.  The \emph{chromatic number} $\chi(X)$ of $X$ is the fewest number of colors necessary to color the vertices of $X$ so that no adjacent vertices are colored the same.  Then a graph $X$ is \emph{perfect} if, for each of its vertex-induced subgraphs $F$, $\omega(F) = \chi(F)$.  Berge conjectured in 1960 that a graph $X$ is perfect if and only if $s_v(C_k, X) = 0$ and $s_v(\bar{C_k}, X) = 0$ for all odd $k > 4$.  The early history of this conjecture can be found in \cite{MR1861355}.  In 1988, Chv{\'a}tal and Sbihi \cite{MR930204} called graphs which satisfied $s_v(C_k, X) = 0$ and $s_v(\bar{C_k}, X) = 0$ for all odd $k > 4$ \emph{Berge graphs.}  Recently, Berge's conjecture was proven by Chudnovsky, Robertson, Seymour, and Thomas \cite{MR2233847}.  Their result is now known as the strong perfect graph theorem and is one of the most important results in recent mathematics.\\

A related class of graphs is those which are chordal.  \emph{Chordal graphs} are those for which $s_v(C_k, X) = 0$ for all $k > 4$.  Chordal graphs have some very interesting properties.  For instance, many problems, such as \emph{minimum coloring, maximum clique, maximum independent set}, and \emph{minimum covering by cliques}, which are NP-complete in general can be solved in polynomial time \cite{MR0327580} on chordal graphs.  In addition, every chordal graph is perfect.\\

Based upon the definitions' reliance upon $s_v(C_k, X) = 0$, it is no surprise that edge zeta functions can distinguish these graph classes:\\
\begin{Cor}[holes, antiholes, perfect, and chordal]
Let $X$ be a graph on $n$ vertices.  Then for $r = 3, \cdots, n$, the numbers $s_v(C_r, X)$ and $s_v(\bar{C_r}, X)$ are determined.  In particular, the edge zeta function can determine whether a graph is chordal or perfect.\\
\end{Cor}

One strong reason for studying edge zeta functions is that they generalize the Ihara zeta function.  We might hope that some of the properties determined so easily by the edge zeta function might, in fact, be determined by the Ihara zeta function.\\

\begin{figure}[t]
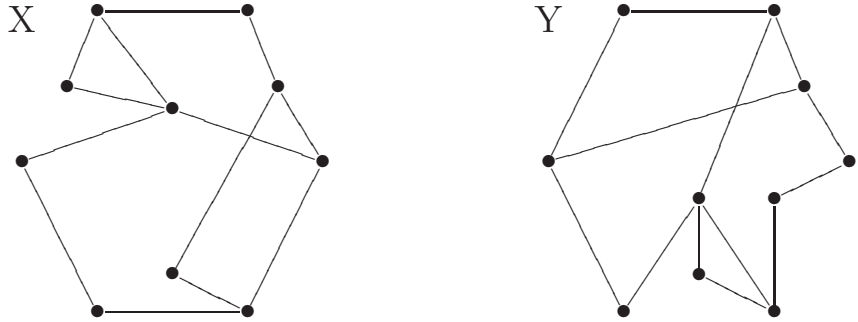

\[\xy
(0,0)*!C++\xybox{
(0, 19)*{\text{\large{X}}};
(0,0)*{\bullet} = "v5";
(30,20)*{\bullet} = "v1";
(6,10)*{\bullet} = "v2";
(40,0)*{\bullet} = "v3";
(30,-20)*{\bullet} = "v8";
(10,-20)*{\bullet} = "v0";
(10,20)*{\bullet} = "v6";
(34,10)*{\bullet} = "v7";
(20,-15)*{\bullet} = "v4";
(20,7)*{\bullet} = "v9";
{\ar@{-} "v0"; "v8"};
{\ar@{-} "v0"; "v5"};
{\ar@{-} "v1"; "v6"};
{\ar@{-} "v1"; "v7"};
{\ar@{-} "v2"; "v9"};
{\ar@{-} "v2"; "v6"};
{\ar@{-} "v3"; "v8"};
{\ar@{-} "v3"; "v9"};
{\ar@{-} "v3"; "v7"};
{\ar@{-} "v4"; "v8"};
{\ar@{-} "v4"; "v7"};
{\ar@{-} "v5"; "v9"};
{\ar@{-} "v6"; "v9"};
};
(70,0)*!C++\xybox{
(0, 19)*{\text{\large{Y}}};
(0,0)*{\bullet} = "v0";
(30,-20)*{\bullet} = "v1";
(30,20)*{\bullet} = "v2";
(40,0)*{\bullet} = "v3";
(10,-20)*{\bullet} = "v4";
(20,-15)*{\bullet} = "v5";
(10,20)*{\bullet} = "v6";
(30,-5)*{\bullet} = "v7";
(34,10)*{\bullet} = "v8";
(20,-5)*{\bullet} = "v9";
{\ar@{-} "v0"; "v8"};
{\ar@{-} "v0"; "v4"};
{\ar@{-} "v0"; "v6"};
{\ar@{-} "v1"; "v9"};
{\ar@{-} "v1"; "v5"};
{\ar@{-} "v1"; "v7"};
{\ar@{-} "v2"; "v8"};
{\ar@{-} "v2"; "v9"};
{\ar@{-} "v2"; "v6"};
{\ar@{-} "v3"; "v8"};
{\ar@{-} "v3"; "v7"};
{\ar@{-} "v4"; "v9"};
{\ar@{-} "v5"; "v9"};
};
\endxy\]
\caption{Two graphs with the same Ihara zeta function.}
\label{fig:SameZeta}
\end{figure}

\begin{Example}[Same zeta function but different structures]
In Figure \ref{fig:SameZeta}, we have an example of two connected md2 graphs which have the same Ihara zeta function.  They both satisfy $\omega(X) = \omega(Y) = 3$.  However, $Y$ is Hamiltonian (with Ham$(Y) = 1$), and $X$ is not.  In addition, we have $s_v(C_6, X) = s_v(C_7, X) = 0$, and $s_v(C_6, Y) = s_v(C_7, Y) = 1$.\\

These graphs were found as part of an effort to enumerate graphs with the same zeta function using McKay's program nauty \cite{NAUTY}.  They were identified, and the data was evaluated, using code written by the author in SAGE \cite{SAGE}.\\
\end{Example}

From this example, we suspect that none of the corollaries in this section are true, in general, for the Ihara zeta function.  We leave it as a problem to find an example of two graphs with the same zeta function where one is perfect or chordal and the other is not or where they have differing clique numbers.\\

We conclude by noting that it is not necessarily a bad thing that the Ihara zeta function does not determine these invariants.  The Ihara zeta function seems to do a decent job at distinguishing md2 graphs, and it may be that it lands in a blind spot that other graph invariants are unable to see.  By combining it with other graph invariants, we have very high hopes for its ability to distinguish graphs.

%\bibliographystyle{plain}
%\bibliography{../AllArticles}

\end{document}